\documentclass[12pt]{article}
\usepackage{amsmath,amsfonts,amssymb,amscd}

\newtheorem{theo}{Theorem}[section]

\newtheorem{lemm}[theo]{Lemma}

\makeatletter \@addtoreset{figure}{section} \makeatother
\makeatletter
\long\def\@makecaption#1#2{%
   \vskip 10\p@
   \setbox\@tempboxa\hbox{{#1}\ \ #2}%
   \ifdim \wd\@tempboxa >\hsize
       {#1}\ \ #2\par
   \else
       \hbox to\hsize{\hfil\box\@tempboxa\hfil}%
   \fi}
\makeatother

\def\qed{\hfill \rule{4pt}{7pt}}
\def\pf{\noindent {\it Proof.} }
\title{ Color degree and color neighborhood union conditions
for long heterochromatic paths in edge-colored graphs \footnote{
Research supported by NSFC. } }
\author{
\small  He Chen and Xueliang Li\\
[2mm]
\small Center for Combinatorics and LPMC \\
\small Nankai University, Tianjin 300071, China \\
\small lxl@nankai.edu.cn\\
}
\date{ }
\begin{document}
\maketitle \vskip5mm
\begin{abstract}
Let $G$ be an edge-colored graph. A heterochromatic (rainbow, or
multicolored) path of $G$ is such a path in which no two edges
have the same color. Let $d^c(v)$ denote the color degree and
$CN(v)$ denote the color neighborhood of a vertex $v$ of $G$. In a
previous paper, we showed that if $d^c(v)\geq k$ (color degree
condition) for every vertex $v$ of $G$, then $G$ has a
heterochromatic path of length at least
$\lceil\frac{k+1}{2}\rceil$, and if $|CN(u)\cup CN(v)|\geq s$
(color neighborhood union condition) for every pair of vertices
$u$ and $v$ of $G$, then $G$ has a heterochromatic path of length
at least $\lceil\frac{s}{3}\rceil+1$. Later, in another paper we
first showed that if $k\leq 7$, $G$ has a heterochromatic path of
length at least $k-1$, and then, based on this we use induction on
$k$ and showed that if $k\geq 8$, then $G$ has a heterochromatic
path of length at least $\lceil\frac{3k}{5}\rceil+1$. In the
present paper, by using a simpler approach we further improve the
result by showing that if $k\geq 8$, $G$ has a heterochromatic
path of length at least $\lceil\frac{2k}{3}\rceil+1$, which
confirms a conjecture by Saito. We also improve a previous result
by showing that under the color neighborhood union condition, $G$
has a heterochromatic path of length at least $\lfloor\frac{2s+4}{5}\rfloor$.\\
[2mm] {\bf Keywords:} edge-colored graph, color degree, color
neighborhood, heterochromatic (rainbow, or multicolored) path.\\
[2mm] {\bf AMS Subject Classification (2000)}: 05C38, 05C15
\end{abstract}

\section{Introduction}

We use Bondy and Murty \cite{B-M} for terminology and notations
not defined here and consider simple graphs only.

Let $G=(V,E)$ be a graph. By an {\it edge-coloring} of $G$ we will
mean a function $C: E\rightarrow \mathbb{N} $, the set of natural
numbers. If $G$ is assigned such a coloring, then we say that $G$
is an {\it edge-colored graph}. Denote the colored graph by
$(G,C)$, and call $C(e)$ the {\it color} of the edge $e\in E$. We
say that $C(uv)=\emptyset$ if $uv\notin E(G)$ for $u,v\in V(G)$.
For a subgraph $H$ of $G$, we denote $C(H)=\{C(e) \ | \ e\in
E(H)\}$ and $c(H)=|C(H)|$. For a vertex $v$ of $G$, the {\it color
neighborhood} $CN(v)$ of $v$ is defined as the set $\{C(e)\ | \ e
\mbox{ is incident with }v\}$ and the {\it color degree} is
$d^c(v)=|CN(v)|$. A path is called {\it heterochromatic (rainbow,
or multicolored)} if any two edges of it have different colors. If
$u$ and $v$ are two vertices on a path $P$, $uPv$ denotes the
segment of $P$ from $u$ to $v$, whereas $vP^{-1}u$ denotes the
same segment but from $v$ to $u$.

There are many existing literature dealing with the existence of
paths and cycles with special properties in edge-colored graphs.
In \cite{C-M-M}, the authors showed that for a 2-edge-colored
graph $G$ and three specified vertices $x, y$ and $z$, to decide
whether there exists a color-alternating path from $x$ to $y$
passing through $z$ is NP-complete. The heterochromatic
Hamiltonian cycle or path problem was studied by Hahn and
Thomassen \cite{H-T}, R\"{o}dl and Winkler (see \cite{F-R}),
Frieze and Reed \cite{F-R}, and Albert, Frieze and Reed
\cite{A-F-R}. For more references, see
\cite{A-J-Z,E-T-1,E-T-2,M-S-T,M-S-T-V}. Many results in these
papers were proved by using probabilistic methods.

In \cite{B-L}, the authors showed that if $G$ is an edge-colored
graph with $d^c(v)\geq k$ (color degree condition) for every
vertex $v$ of $G$, then $G$ has a heterochromatic path of length
at least $\lceil\frac{k+1}{2}\rceil$, and if $|CN(u)\cup
CN(v)|\geq s$ (color neighborhood union condition) for every pair
of vertices $u$ and $v$ of $G$, then $G$ has a heterochromatic
path of length at least $\lceil\frac{s}{3}\rceil+1$. In
\cite{C-L}, we first showed that if $3\leq k\leq 7$, $G$ has a
heterochromatic path of length at least $k-1$, and then, based on
this we use induction on $k$ and showed that if $k\geq 8$, then
$G$ has a heterochromatic path of length at least
$\lceil\frac{3k}{5}\rceil+1$. In the present paper, by using a
simpler approach we further improve the result by showing that if
$k\geq 8$, $G$ has a heterochromatic path of length at least
$\lceil\frac{2k}{3}\rceil+1$, which confirms a conjecture by
Saito. We also show that under the color neighborhood union
condition, $G$ has a heterochromatic path of length at least
$\lfloor\frac{2s+4}{5}\rfloor$.

\section{Long heterochromatic paths under the color degree condition}

In this section we will give a better lower bound for the length
of the longest heterochromatic path in $G$ when $k\geq 8$. As an
induction initial, we need the following result as a lemma.

\begin{lemm} [\cite{C-L}]  Let $G$ be an edge-colored graph and
$3\leq k\leq 7$ an integer. Suppose that $d^c(v)\geq k$ for every
vertex $v$ of $G$. Then $G$ has a heterochromatic path of length
at least $k-1$.
\end{lemm}

Then, we need to do the following preparations.
\begin{lemm}Suppose $P=u_1u_2\ldots u_lu_{l+1}$ is a longest heterochromatic
path. If there exists an $x$ such that $3\leq x\leq l$ and
$C(u_1u_x)\notin C(P)$, then
$C(u_{x-1}u_x)\notin(CN(u_{l+1})-C(u_1u_{l+1},u_2u_{l+1},\ldots,u_{l-1}u_{l+1}))$.
\end{lemm}
{\bf\pf}By contradiction. If there exists a $v\in V(G)-V(P)$ such
that $C(u_{l+1}v)=C(u_{x-1}u_x)$, then
$u_{x-1}P^{-1}u_1u_xPu_{l+1}v$ is a heterochromatic path of length
$l+1$, a contradiction. So
$C(u_{x-1}u_x)\notin(CN(u_{l+1})-C(u_1u_{l+1},u_2u_{l+1},\ldots,u_{l-1}u_{l+1}))$.
\qed
\begin{lemm}Suppose $P=u_1u_2\ldots u_lu_{l+1}$ is a longest
heterochromatic path, $v\in V(G)-V(P)$ and
$C(u_{l+1}v)=C(u_1u_2)$. If there exists an $x$ such that $2\leq
x\leq l-2$ and $|C(u_xv,u_{x+2}v)-C(P)|=2$, then
$C(u_xu_{x+1},u_{x+1}u_{x+2})\cap(CN(u_{l+1})-C(u_1u_{l+1},u_2u_{l+1},
\ldots,u_{l-1}u_{l+1}))=\emptyset$.
\end{lemm}
{\bf\pf}By contradiction. If there exists a $v'\in V(G)-V(P)$ such
that $u_{l+1}v'\in E(G)$ and $C(u_{l+1}v')\in
C(u_xu_{x+1},u_{x+1}u_{x+2})$, then $u_1Pu_xvu_{x+2}Pu_{l+1}v'$ is
a heterochromatic path of length $l+1$, a contradiction. So
$C(u_xu_{x+1},u_{x+1}u_{x+2})\cap(CN(u_{l+1})-C(u_1u_{l+1},u_2u_{l+1},
\ldots,u_{l-1}u_{l+1}))=\emptyset$.
\qed
\begin{lemm}Let $P=u_1u_2\ldots u_lu_{l+1}v_1$ be a path in $G$ such
that\\
\indent (a) $u_1Pu_{l+1}$ is a longest heterochromatic path in $G$;\\
\indent (b) $C(u_{l+1}v_1)=C(u_{j_0}u_{j_0+1})$ and $1\leq j_0\leq
l$ is as small as possible, subject to (a).\\
Then $C(u_1u_{j_0+1},u_1u_{j_0+2},\ldots,u_1u_{2j_0})\subseteq
C(P)$.
\end{lemm}
{\bf\pf}By contradiction. If there exists an $x$ such that
$j_0+1\leq x\leq 2j_0$ and $C(u_1u_x)\notin C(P)$, then
$u_{x-1}P^{-1}u_1u_xPu_{l+1}$ is a heterochromatic path of length
$l$ and $u_{j_0+1}u_{j_0}$ is the $x-j_0-1\leq
2j_0-j_0-1=j_0-1<j_0$-th edge in this heterochromatic path,
contradicting the choice of $P$. Therefore
$C(u_1u_{j_0+1},u_1u_{j_0+2},\ldots,u_1u_{2j_0})\subseteq C(P)$.
\qed
\begin{lemm} Let $P=u_1u_2\ldots u_lu_{l+1}v_1$ be a path in $G$ such
that\\
\indent (a) $u_1Pu_{l+1}$ is a longest heterochromatic path in $G$;\\
\indent (b) $C(u_{l+1}v_1)=C(u_{j_0}u_{j_0+1})$ and $1\leq j_0\leq
l$ is as small as possible, subject to (a).\\
Then for any $2j_0+1\leq x\leq l$,
$|C(u_1u_x,u_1u_{x+1})-C(P)|\leq 1$.
\end{lemm}
{\bf\pf}By induction. If there exists an $x$ such that $2j_0+1\leq
x\leq l$ and $|C(u_1u_x,u_1u_{x+1})-C(P)|=2$, then
$u_2Pu_xu_1u_{x+1}Pu_{l+1}$ is a heterochromatic path of length
$l$ and $u_{j_0}u_{j_0+1}$ is the $(j_0-1)$-th edge in this
heterochromatic path, contradicting the choice of $P$. Therefore
$|C(u_1u_x,u_1u_{x+1})-C(P)|\leq 1$ for any $2j_0+1\leq x\leq
l$.\qed
\begin{lemm}Suppose $d^c(v)\geq k$ for every vertex $v\in V(G)$ and
the length of a longest heterochromatic path in $G$ is
$l=\lceil\frac{2k}{3}\rceil$. Then there is a heterochromatic path
$P=u_1u_2\ldots u_lu_{l+1}$ in $G$ and a $v\in V(G)-V(P)$ such
that $C(u_{l+1}v)=C(u_1u_2)$.
\end{lemm}
{\bf\pf}Let $P=u_1u_2\ldots u_lu_{l+1}v_1$ be a path in $G$ such
that\\
\indent (a) $u_1Pu_{l+1}$ is a longest heterochromatic path in $G$;\\
\indent (b) $C(u_{l+1}v_1)=C(u_{j_0}u_{j_0+1})$ and $1\leq j_0\leq
l$ is as small as possible, subject to (a).\\
\indent Then we claim that $j_0=1$. We will show this by
contradiction. Suppose $j_0>1$. Denote $i_j=C(u_ju_{j+1})$ for
$1\leq j\leq l$.

Since the longest heterochromatic path in $G$ is of length $l$,
for any $v\in V(G)-V(u_1Pu_{l+1})$ we have $C(u_1v)\in C(P)$. On
the other hand, $d^c(u_1)\geq k$. So there are at least
$k-l=\lfloor\frac{k}{3}\rfloor$ different colors not in $C(P)$
appearing in $\{u_1u_3,u_1u_4,\ldots,u_1u_l,u_1u_{l+1}\}$. Then
there are $x_i$'s such that $3\leq x_1<x_2<\ldots<x_{k-l}\leq l+1$
and
$|C(\{u_1u_{x_1},u_1u_{x_2},\ldots,u_1u_{x_{k-l}}\})-C(P)|=k-l$.
Therefore, by Lemma 2.2 and the assumption that $j_0>1$ we have
$(CN(u_{l+1})-C(u_1u_{l+1},u_2u_{l+1},\ldots,u_{l-1}u_{l+1},u_lu_{l+1}))\subseteq
C(P)-\{i_1,i_{x_1-1},\ldots,i_{x_{k-l}-1}\}$. Since
$d^c(u_{l+1})\geq k$, we have $\lceil\frac{2k}{3}\rceil=l\geq
|C(u_1u_{l+1},u_2u_{l+1},\ldots,u_lu_{l+1})|\geq
k-|C(P)-\{i_1,i_{x_1-1},\ldots,i_{x_{k-l}-1}\}| = k-(l-k+l-1) =
2k-2l+1 =2\lfloor\frac{k}{3}\rfloor+1$. Since if $k\equiv 0 \ (mod
 \ 3)$ then $2\lfloor\frac{k}{3}\rfloor+1>\lceil\frac{2k}{3}\rceil$, we
need only to consider the cases when $k\equiv 1 \ (mod \ 3)$ or
$k\equiv 2 \ (mod \ 3)$.

{\bf Case 1.} $k\equiv 1 \ (mod \ 3)$.\\
\indent In this case, we have
$\lceil\frac{2k}{3}\rceil=2\lfloor\frac{k}{3}\rfloor+1$. Then
$CN(u_{l+1})-C(u_1u_{l+1},u_2u_{l+1},\ldots,$
$u_lu_{l+1})=C(P)-\{i_1,i_{x_1-1},\ldots,i_{x_{k-l}-1}\}$ and
$C(u_lu_{l+1})\notin C(P)-\{i_1,i_{x_1-1},\ldots,$
$i_{x_{k-l}-1}\}$. Then, we can get
$C(u_lu_{l+1})=i_{x_{k-l}-1}$, i.e., $x_{k-l}=l+1$.\\
\indent Noticing that $C(u_{l+1}v_1)=i_{j_0}$ and $j_0$ is as
small as possible, we have $\{i_2,\ldots,$
$i_{j_0-1}\}\cap(C(P)-\{i_1,i_{x_1-1},\ldots,i_{x_{k-l}-1}\})=\emptyset$,
and so $\{3,4,\ldots,j_0\}\subseteq\{x_1,x_2,\ldots,$ $x_{k-l}\}$.
Hence, by Lemmas 3.3 and 3.4 we have
$\{x_1,x_2,\ldots,x_{k-l}\}\subseteq\{3,4,\ldots,j_0\}\cup\{2j_0+1,\ldots,l+1\},$
and
$|\{x_1,x_2,\ldots,x_{k-l}\}\cap\{2j_0+1,2j_0+2,\ldots,l+1\}|\leq
\lfloor\frac{(l+1)-(2j_0+1)}{2}\rfloor+1=\lfloor
\frac{l}{2}\rfloor-j_0+1$. Consequently,
$|\{x_1,x_2,\ldots,x_{k-l}\}|\leq (j_0-2)+\lfloor
\frac{l}{2}\rfloor-j_0+1=\lfloor\frac{l}{2}\rfloor-1<\lfloor\frac{k}{3}\rfloor=k-l$,
a contradiction.

{\bf Case 2.} $k\equiv 2 \ (mod \ 3)$.\\
\indent In this case, we have
$\lceil\frac{2k}{3}\rceil=(2\lfloor\frac{k}{3}\rfloor+1)+1$. We
distinguish the following two cases:

{\bf Case 2.1.} $x_{k-l}=l+1$.\\
\indent Since $\{x_1,x_2,\ldots,x_{k-l}\}\subseteq
\{3,4,\ldots,j_0\}\cup\{2j_0+1,\ldots,l+1\}$ by Lemma 2.4, and
$|\{x_1,x_2,\ldots,x_{k-l}\}\cap\{3,4,\ldots,j_0\}|\leq
|\{3,4,\ldots,j_0\}|=j_0-2$,
$|\{x_1,x_2,\ldots,x_{k-l}\}\cap\{2j_0+1,\ldots,l+1\}|\leq
\lfloor\frac{(l+1)-(2j_0+1)}{2}\rfloor+1=\lfloor\frac{l}{2}\rfloor-j_0+1$
by Lemma 2.5, we have $|\{x_1,x_2,\ldots,x_{k-l}\}|\leq
(j_0-2)+(\lfloor\frac{l}{2}\rfloor-j_0+1)=\lceil\frac{l}{2}\rceil-1=\lfloor\frac{k}{3}\rfloor=k-l$.
Then
$\{x_1,x_2,\ldots,x_{k-l}\}=\{3,4,\ldots,j_0,2j_0+1,2j_0+3,\ldots,l-1,l+1\}$.\\
\indent Since
$\lceil\frac{2k}{3}\rceil=(2\lfloor\frac{k}{3}\rfloor+1)+1$, there
is at most one color in $C(P)-\{i_1,i_{x_1-1},\ldots,$
$i_{x_{k-l}-1}\}$ contained in
$C(\{u_1u_{l+1},u_2u_{l+1},\ldots,u_lu_{l+1}\})$, i.e.,
$|(C(P)-\{i_1,i_{x_1-1},\ldots,$
$i_{x_{k-l}-1}\})-(CN(u_{l+1})-C(\{u_1u_{l+1},u_2u_{l+1},
\ldots,u_lu_{l+1}\}))|\leq 1$.\\
\indent If $j_0\geq 3$, then $|\{j_0+1,\ldots,2j_0-1\}|=j_0-1\geq
2$, and so there exists a $v\in V(G)-V(P)$ such that
$C(u_{l+1}v)=i_{j_0+s}$ for some $1\leq s\leq j_0-1$. Then
$u_{2j_0}P^{-1}u_1u_{2j_0+1}Pu_{l+1}$ is a heterochromatic path of
length $l$ and $u_{j_0+s+1}u_{j_0+s}$ is the $(j_0-s)$-th edge in
this heterochromatic path, contradicting the choice of $P$.\\
\indent Therefore we need only to consider the case when $j_0=2$,
then $x_1=2j_0+1=5$. In this case, there exists a $v\in V(G)-V(P)$
such that $C(u_{l+1}v)\in\{i_{j_0+1},i_{x_1}\}$. If
$C(u_{l+1}v)=i_{j_0+1}=i_{2j_0-1}$, then $u_{2j_0}u_{2j_0-1}$ is
the first edge in the heterochromatic path
$u_{2j_0}P^{-1}u_1u_{2j_0+1}Pu_{l+1}$ of length $l$; if
$C(u_{l+1}v)=i_{x_1}$, then $u_{x_1+1}u_{x_1}$ is the first edge
in the heterochromatic path $u_{x_1+1}P^{-1}u_1u_{x_2}Pu_{l+1}$ of
length $l$, contradicting the choice of $P$.

{\bf Case 2.2.} $x_{k-l}<l+1$.\\
\indent In this case, we can get
$\{x_1,x_2,\ldots,x_{k-l}\}\subseteq\{3,4,\ldots,j_0\}\cup\{2j_0+1,2j_0+2,\ldots,l\}$
by Lemma 3.3,
$|\{x_1,x_2,\ldots,x_{k-l}\}\cap\{3,4,\ldots,j_0\}|\leq j_0-2$ and
$|\{x_1,x_2,\ldots,$ $x_{k-l}\}\cap\{2j_0+1,\ldots,l\}|\leq
\lfloor\frac{l-(2j_0+1)}{2}\rfloor+1=\lfloor\frac{l-1}{2}\rfloor-j_0+1=\frac{l}{2}-j_0$
by Lemma 3.4. Consequently, $|\{x_1,x_2,\ldots,x_{k-l}\}|\leq
(j_0-2)+(\frac{l}{2}-j_0)=\frac{l}{2}-2=\lfloor\frac{k}{3}\rfloor-1<k-l$,
a contradiction.

From the arguments of all the above cases, we get that $j_0$
cannot be larger than $1$, and so $j_0=1$. \qed

Now we are ready to give our main result.

\begin{theo}If $d^c(v)\geq k\geq 7$ for any $v\in V(G)$, then
$G$ has a heterochromatic path of length at least
$\lceil\frac{2k}{3}\rceil+1$.
\end{theo}
{\bf\pf}We will prove the theorem by induction.

If $k=7$, our Lemma 2.1 guarantees that $G$ has a heterochromatic
path of length at least $6=\lceil\frac{2\times 7}{3}\rceil+1$.

Assume that if $d^c(v)\geq k-1$ for any $v\in V(G)$, $G$ has a
heterochromatic path of length at least
$\lceil\frac{2(k-1)}{3}\rceil+1$. Then we need only to show that
if $d^c(v)\geq k$ for any $v\in V(G)$, $G$ has a heterochromatic
path of length $\lceil\frac{2k}{3}\rceil+1$. Since if $k\equiv 0 \
(mod  \ 3)$ then
$\lceil\frac{2(k-1)}{3}\rceil+1=\lceil\frac{2k}{3}\rceil+1$, we
need only to show that if $k\equiv 1,2 \ (mod \ 3)$, $G$ has a
heterochromatic path of length at least
$\lceil\frac{2k}{3}\rceil+1$.

By the assumption we know that $G$ has a heterochromatic path of
length at least
$\lceil\frac{2(k-1)}{3}\rceil+1=\lceil\frac{2k}{3}\rceil$. Assume
that the longest heterochromatic path in $G$ is of length
$\lceil\frac{2k}{3}\rceil$. Then, by Lemma 2.6 $G$ has a
heterochromatic path $P=u_1u_2\ldots u_lu_{l+1}$ of length
$l=\lceil\frac{2k}{3}\rceil$ and there exists a $v_1\in V(G)-V(P)$
such that $C(u_{l+1}v_1)=C(u_1u_2)$. Denote $i_j=C(u_ju_{j+1})$
for $1\leq j\leq l$.

Since $d^c(v_1)\geq k$, we have that $d^c(u_1)\geq k$ and the
longest heterochromatic path in $G$ is of length $l$, and so there
exist $y_i$'s and $x_j$'s such that $2\leq
y_1<y_2<y_3<\ldots<y_{k-l}\leq l$ and $3\leq
x_1<x_2<\ldots<x_{k-l}\leq l+1$, and
$|C(\{u_{y_1}v_1,u_{y_2}v_1,\ldots,u_{y_{k-l}}v_1\})-C(P)|=k-l$
and
$|C(\{u_1u_{x_1},u_1u_{x_2}v_1,\ldots,u_1u_{x_{k-l}}\})$ $-C(P)|=k-l$.\\
\indent If there exists a $j_0$ such that $1\leq j_0\leq k-l-1$
and $y_{j_0+1}=y_{j_0}+1$, then
$u_1Pu_{y_{j_0}}v_1u_{y_{j_0}+1}Pu_{l+1}$
is a heterochromatic path of length $l+1$, a contradiction.\\
\indent If $y_1=2$, since $k-l=\lfloor\frac{k}{3}\rfloor\geq 2$,
then there exists a $j$ such that $1\leq j\leq k-l$ and
$C(u_1u_{x_j})\notin C(P)\cup C(u_2v_1)$. Then
$u_{x_j-1}P^{-1}u_2v_1u_{l+1}P^{-1}u_{x_j}u_1$ is a
heterochromatic path of length $l+1$, a contradiction.\\
\indent If $x_{k-l}=l+1$, since $k-l=\lfloor\frac{k}{3}\rfloor\geq
2$, then there exists a $j$ such that $1\leq j\leq k-l$ and
$C(u_{y_j}v_1)\notin(C(P)\cup C(u_1u_{l+1}))$. Then
$v_1u_{y_j}P^{-1}u_1u_{l+1}P^{-1}u_{y_j+1}$ is a heterochromatic
path of length $l+1$, a contradiction.\\
\indent If there exists a $j_0$ such that $1\leq j_0\leq k-l-1$
and $x_{j_0+1}=x_{j_0}+1$, then
$u_2Pu_{x_{j_0}}u_1u_{x_{j_0}+1}Pu_{l+1}v_1$ is a
heterochromatic path of length $l+1$, a contradiction.\\
\indent Consequently, we have $3\leq x_1<x_1+1<x_2<x_2+1<\ldots
<x_{k-l}\leq l$ and $3\leq y_1<y_1+1<y_2<y_2+1<\ldots<y_{k-l}\leq
l$. Then
$CN(u_{l+1})-C(\{u_1u_{l+1},u_2u_{l+1},\ldots,u_{l-1}u_{l+1},
u_lu_{l+1}\})\subseteq
C(P)-\{i_{x_1-1},i_{x_2-1},\ldots,i_{x_{k-l}-1}\}$ by Lemma 2.2
and the fact that $P$ is the longest heterochromatic path in $G$.
On the other hand, $x_{k-l}\leq l$, and so $i_l\in
C(P)-\{i_{x_1-1},i_{x_2-1},\ldots,i_{x_{k-l}-1}\}$. Then
$CN(u_{l+1})-C(\{u_1u_{l+1},u_2u_{l+1},\ldots,u_{l-1}u_{l+1}\})\subseteq
C(P)-\{i_{x_1-1},i_{x_2-1},\ldots,i_{x_{k-l}-1}\}$.

We distinguish the following two cases:

{\bf Case 1.} $k\equiv 1 \ (mod \ 3)$.\\
\indent Since $d^c(u_{l+1})\geq k$, we have $l-1\geq
|C(\{u_1u_{l+1},u_2u_{l+1},\ldots,u_{l-1}u_{l+1}\})|\geq
k-|C(P)-\{i_{x_1-1},i_{x_2-1},\ldots,i_{x_{k-l}-1}\}|=k-l+(k-l) =
l-1$. Therefore $u_1u_{l+1}\in E(G)$ and $C(u_1u_{l+1})\in
\{i_{x_1-1},i_{x_2-1},\ldots,i_{x_{k-l}-1}\}$. Suppose
$C(u_1u_{l+1})=C(u_{x_j-1}u_{x_j})$ for some $1\leq j\leq k-l$.\\
\indent On the other hand, since $3\leq
y_1<y_1+1<y_2<y_2+1<\ldots<y_{k-l}\leq l$ and $3\leq
x_1<x_1+1<x_2<x_2+1<\ldots<x_{k-l}\leq l$, we get that
$2(k-l)-2\leq y_{k-l}-y_1\leq l-3=2(k-l)-2$ and $2(k-l)-2\leq
x_{k-l}-x_1\leq l-3=2(k-l)-2$, and then
$\{y_1,y_2,\ldots,y_{k-l}\}=\{3,5,\ldots,l-2,l\}=\{x_1,x_2,\ldots,x_{k-l}\}$,
$v_1u_{x_j}Pu_{l+1}u_1Pu_{x_j-1}$ is a heterochromatic path of
length $l+1$, a contradiction.

{\bf Case 2.} $k\equiv 2 \ (mod \ 3)$.\\
\indent Since $3\leq y_1<y_1+1<y_2<y_2+1<\ldots<y_{k-l}\leq l$, we
have $2(k-l-1)\leq y_{k-l}-y_1\leq l-3 = 2(k-l-1)+1$. Then we get
that $y_{j+1}=y_j+2$ for $j=1,2,\ldots,k-l-1$ or there exists a
$j_0$ such that $1\leq j_0\leq k-l-1$, and $y_{j+1}=y_j+2$ for any
$1\leq j\leq k-l-1$ and $j\neq j_0$, $y_{j_0+1}=y_{j_0}+3$.

{\bf Case 2.1} $y_{j+1}=y_j+2$ for $j=1,2,\ldots, k-l-1$.\\
\indent In this case, we have
$(CN(u_{l+1})-C(\{u_1u_{l+1},u_2u_{l+1},\ldots,u_{l-1}u_{l+1},u_lu_{l+1}\}))
\subseteq C(P)-\{i_{y_1},i_{y_1+1}, i_{y_2},
i_{y_2+1},\ldots,i_{y_{k-l}-1}\}$ by Lemma 2.3 and the fact that
$P$ is the longest heterochromatic path in $G$. Noticing that
$y_{k-l}\leq l$, we have $i_l\notin \{i_{y_1},i_{y_1+1}, i_{y_2},
i_{y_2+1},\ldots,i_{y_{k-l}-1}\}$. Then
$(CN(u_{l+1})-C(\{u_1u_{l+1},u_2u_{l+1},\ldots, u_{l-1}u_{l+1}\\
\})) \subseteq C(P)-\{i_{y_1},i_{y_1+1}, i_{y_2},
i_{y_2+1},\ldots,i_{y_{k-l}-1}\}$.\\
\indent On the other hand,
$CN(u_{l+1})-C(\{u_1u_{l+1},u_2u_{l+1},\ldots,u_{l-1}u_{l+1}\})\subseteq
C(P)-\{i_{x_1-1}, i_{x_2-1},\ldots,i_{x_{k-l}-1}\}$. Therefore
$(CN(u_{l+1})-C(\{u_1u_{l+1},u_2u_{l+1},\ldots,u_{l-1}u_{l+1}\}\\
)) \subseteq C(P)-\{i_{y_1},i_{y_1+1}, i_{y_2},
i_{y_2+1},\ldots,i_{y_{k-l}-1}\}\cup\{i_{x_1-1},i_{x_2-1},
\ldots,i_{x_{k-l}-1}\}$.\\
\indent Note that $3\leq x_1<x_1+1<x_2<\ldots<x_{k-l}\leq l$. Then
$\{x_1,x_2,\ldots,x_{k-l}\}-\{y_1+1,y_2,y_2+1,\ldots,y_{k-l}\}\neq\emptyset$
and
$\{i_{x_1-1},i_{x_2-1},\ldots,i_{x_{k-l}-1}\}-\{i_{y_1},i_{y_1+1},i_{y_2},
\ldots,\\ i_{y_{k-l}-1}\}\neq\emptyset$.\\
\indent Consequently, $l-1\geq
|C(\{u_1u_{l+1},u_2u_{l+1},\ldots,u_{l-1}u_{l+1}\})|\geq
k-|C(P)-\{i_{y_1},i_{y_1+1}, i_{y_2},
i_{y_2+1},\ldots,i_{y_{k-l}-1}\}\cup\{i_{x_1-1},i_{x_2-1},
\ldots,i_{x_{k-l}-1}\}|\geq k-l+2(k-l-1)+1 = 3k-3l-1$. It is easy
to check that if $k>8$, $3k-3l-1>l-1$, and so we need only to
consider the case when $k=8$.

If $k=8$, $l-1=3k-3l-1$, and so we need only to consider the case
when $|\{i_{x_1-1},i_{x_2-1}\}-\{i_{y_1},i_{y_1+1}\}|=1$. Denote
$i_7=C(u_{y_1}v_1)$, $i_8=C(u_{y_2}v_1)$. We
distinguish the following two cases:\\
\indent {\bf Case 2.1.1} $y_1=3$ and $y_2=5$.\\
\indent In this case, we need only to consider the cases when
$x_1=3$ and $x_2=5$, or $x_1=4$ and $x_2=6$. Then $C(u_1u_7)\in
\{i_2,i_3,i_4\}$ or $C(u_1u_7)\in\{i_3,i_4,i_5\}$. If
$C(u_1u_7)=i_3$ or $i_5$, then $u_4u_5v_1u_3u_2u_1u_7u_6$ is a
heterochromatic path of length $7$; if $C(u_1u_7)=i_2$ or $i_4$,
then $u_4u_3v_1u_5u_6u_7u_1u_2$ is a heterochromatic path of
length $7$, a contradiction.\\
\indent {\bf Case 2.1.2} $y_1=4$ and $y_2=6$.\\
\indent In this case, we need only to consider the cases when
$x_1=3$ and $x_2=5$, or $x_1=3$ and $x_2=6$, or $x_1=4$ $x_2=6$.
Then $C(u_1u_7)\in\{i_2,i_4,i_5\}$ or $\{i_3,i_4,i_5\}$. If
$C(u_1u_7)=i_3$ or $i_5$, then $u_5u_4v_1u_6u_7u_1u_2u_3$ is a
heterochromatic path of length $7$; if $C(u_1u_7)=i_4$, then
$u_5u_6v_1u_4u_3u_2u_1u_7$ is a heterochromatic path of length
$7$. So, we may assume $C(u_1u_7)=i_2$. Then
$C(u_2u_7,u_3u_7,u_4u_7,u_5u_7)\cap\{i_1,i_3,i_4,i_5,i_6\}\subseteq\{i_4,i_5\}$
and $|\{C(u_2u_7,u_3u_7,u_4u_7,u_5u_7)\}|=4$. So $C(u_3u_7)=i_4$
or $i_5$ or some color $\notin\{i_1,i_2,\ldots,i_6\}$. Let
\[
P'=\left\{
\begin{array}{ll}
v_1u_4u_5u_6u_7u_3u_2u_1 &\mbox{if $C(u_3u_7)\notin\{i_1,i_2,\ldots,i_6,i_7\}$;}\\
u_2u_1u_7u_3u_4u_5u_6v_1 &\mbox{if $C(u_3u_7)=i_7$ ;}\\
u_2u_1u_7u_3u_4v_1u_6u_5 &\mbox{if $C(u_3u_7)=i_4$;}\\
u_5u_4v_1u_6u_7u_3u_2u_1 &\mbox{if $C(u_3u_7)=i_5$.}
\end{array}
\right.
\]
Then, $P'$ is a heterochromatic path of length $7$, a
contradiction.

{\bf Case 2.2} There exists a $j_0$ such that $1\leq j_0\leq
k-l-1$, and $y_{j+1}=y_j+2$ for any $1\leq j\leq k-l-1$ and $j\neq
j_0$, $y_{j_0+1}=y_{j_0}+3$.\\
\indent In this case, we have
$CN(u_{l+1})-C(\{u_1u_{l+1},u_2u_{l+1},\ldots,u_{l-1}u_{l+1},u_lu_{l+1}\})\subseteq
C(P)-\{i_{y_1},i_{y_1+1},
\ldots,i_{y_{j_0}-1},i_{y_{j_0+1}},i_{y_{j_0+1}+1},\ldots,i_{y_{k-l}-1}\}$
by Lemma 3.2 and the fact that $P$ is the longest heterochromatic
path in $G$. Note that $y_{k-l}\leq l$, and so $i_l\notin
\{i_{y_1},i_{y_1+1}, i_{y_2}, i_{y_2+1},\ldots,i_{y_{k-l}-1}\}$.
Then
$CN(u_{l+1})-C(\{u_1u_{l+1},u_2u_{l+1},\ldots,u_{l-1}u_{l+1}\}\\
)\subseteq C(P)-\{i_{y_1},i_{y_1+1},
\ldots,i_{y_{j_0}-1},i_{y_{j_0+1}},i_{y_{j_0+1}+1},\ldots,i_{y_{k-l}-1}\}$.\\
\indent On the other hand,
$CN(u_{l+1})-C(\{u_1u_{l+1},u_2u_{l+1},\ldots,u_{l-1}u_{l+1}\})\subseteq
C(P)-\{i_{x_1-1},i_{x_2-1},\ldots,i_{x_{k-l}-1}\}$. Therefore
$(CN(u_{l+1})-C(\{u_1u_{l+1},u_2u_{l+1},\ldots,u_{l-1}u_{l+1}\}))\\
\subseteq C(P)-\{i_{y_1},i_{y_1+1},
\ldots,i_{y_{j_0}-1},i_{y_{j_0+1}},i_{y_{j_0+1}+1},\ldots,i_{y_{k-l}-1}\}\cup\{i_{x_1-1},i_{x_2-1},
\ldots,i_{x_{k-l}-1}\\ \}$.\\
\indent Note that $3\leq x_1<x_1+1<x_2<\ldots<x_{k-l}\leq l$,
$|\{x_1,x_2,\ldots,x_{k-l}\}-\{y_1+1,y_2,y_2+1,\ldots,y_{j_0}\}\cup\{y_{j_0+1}+1,y_{j_0+2},\ldots,y_{k-l}\}|\geq
2$. So $|\{i_{x_1-1},i_{x_2-1},\ldots,i_{x_{k-l}-1}\}\\
-\{i_{y_1},i_{y_1+1},\ldots,
i_{y_{j_0}-1},i_{y_{j_0+1}},i_{y_{j_0+1}+1},\ldots,i_{y_{k-l}-1}\}|\geq
2$.\\
\indent Consequently, $l-1\geq
|C(\{u_1u_{l+1},u_2u_{l+1},\ldots,u_{l-1}u_{l+1}\})|\geq
k-|C(P)-\{i_{x_1-1},i_{x_2-1},
\ldots,i_{x_{k-l}-1}\}\cup\{i_{y_1},i_{y_1+1},\ldots,i_{y_{j_0}-1},i_{y_{j_0+1}},
i_{y_{j_0+1}+1},\ldots,i_{y_{k-l}-1}\}| \geq
k-l+2(j_0-1)+2(k-l-j_0-1)+2=3k-3l-2$. It is easy to check that if
$k>11$, $3k-3l-2>l-1$, and so we need only to consider the cases
when $k=8$ or $k=11$.

\indent {\bf Case 2.2.1} $k=8$. In this case, $y_1=3$ and $y_2=6$.
Denote $i_7=C(u_3v_1)$ and $i_8=C(u_6v_1)$. We distinguish
the following cases:\\
\indent {\bf Case 2.2.1.1} $x_1=3$ and $x_2=5$. Then
$|C(u_1u_7,u_2u_7,u_3u_7,u_4u_7,u_5u_7)-\{i_1,i_3,i_5,
i_6\}|\geq 4$.\\
\indent If $C(u_1u_5)\notin\{i_7,i_8\}$, then
$u_4u_5u_1u_2u_3v_1u_6u_7$ is a heterochromatic path of length
$7$, a contradiction. So we may assume
$C(u_1u_5)\in\{i_7,i_8\}$.\\
\indent If $C(u_1u_7)=i_2$, then $v_1u_3u_4u_5u_6u_7u_1u_2$ is a
heterochromatic path of length $7$, a contradiction. \\
\indent If $C(u_1u_7)=i_4$, then since
$|C(u_1u_7,u_2u_7,u_3u_7,u_4u_7,u5u_7)-\{i_1,i_3,i_5,i_6\}|\geq
4$, we have
$C(u_2u_7,u_3u_7,u_4u_7,u_5u_7)-\{i_1,i_3,i_4,i_5,i_6,i_7,i_8\}\neq\emptyset$.
Let
\[
P'=\left\{
\begin{array}{ll}
v_1u_3u_4u_5u_6u_7u_2u_1 &\mbox{if $C(u_2u_7)\notin\{i_1,i_3,i_4,i_5,i_6,i_7,i_8\}$;}\\
u_5u_6v_1u_3u_4u_7u_1u_2 &\mbox{if $C(u_4u_7)\notin\{i_1,i_3,i_4,i_5,i_6,i_7,i_8\}$;}\\
u_4u_3v_1u_6u_5u_7u_1u_2 &\mbox{if
$C(u_5u_7)\notin\{i_1,i_3,i_4,i_5,i_6,i_7,i_8\}$.}
\end{array}
\right.
\]
Then, $P'$ is a heterochromatic path of length $7$, a
contradiction. So $C(u_3u_7)-\{i_1,i_3,i_4,i_5,i_6,i_7,i_8\}\neq
\emptyset$. In this case, $u_2u_1u_5u_4u_3u_7u_6v_1$ is a
heterochromatic path of length $7$ if $C(u_1u_5)=i_7$. Now it
remains to consider the case when $C(u_1u_5)=i_8$. Since
$u_1u_2u_3v_1u_6u_5u_4$ is a heterochromatic path of length $6$,
$C(u_1u_3,u_1u_4,u_1u_6,u_1v_1)-\{i_1,i_2,i_4,i_5,i_6,i_7,i_8\}\neq
\emptyset$. Let
\[
P''=\left\{
\begin{array}{ll}
u_5u_4u_1u_2u_3v_1u_6u_7 &\mbox{if $C(u_1u_4)\notin\{i_1,i_2,i_4,i_5,i_6,i_7,i_8\}$;}\\
u_2u_3v_1u_7u_6u_1u_5u_4 &\mbox{if $C(u_1u_6)\notin\{i_1,i_2,i_4,i_5,i_6,i_7,i_8\}$;}\\
v_1u_1u_2u_3u_4u_5u_6u_7 &\mbox{if $C(u_1v_1)\notin\{i_1,i_2,i_3,i_4,i_5,i_6,i_7,i_8\}$;}\\
u_2u_3u_1v_1u_7u_6u_5u_4 &\mbox{if $C(u_1v_1)=i_3$.}
\end{array}
\right.
\]
Then, $P''$ is a heterochromatic path of length $7$, and so
$C(u_1u_3)\notin\{i_1,i_2,i_4,i_5,i_6,i_7,\\ i_8\}$, i.e.,
$C(u_1u_3)\notin\{i_1,i_2,i_3,i_4,i_5,i_6,i_7,i_8\}$. Denote
$i_9=C(u_1u_3)$. Since $C(u_1v_1,\\ u_1u_4,
u_1u_6)\subseteq\{i_1,i_2,\ldots, i_6,i_8,i_9\}$ and $d^c(u_1)\geq
8$, there exists a $v_2\notin\{u_1,u_2,\ldots,u_7,\\ v_1\}$ such
that $C(u_1v_2)=i_3$. Then, $v_2u_1u_2u_3v_1u_6u_5u_4$ is a
heterochromatic path of length $7$, a contradiction.\\
\indent If $C(u_1u_7)\neq i_4$, then
$|C(u_2u_7,u_3u_7,u_4u_7,u_5u_7)-\{i_1,i_3,i_5,i_6\}|=4$. Note
that if $C(u_1u_7)=i_3$, then $u_4u_5u_6v_1u_3u_2u_1u_7$ is a
heterochromatic path of length $7$; if $C(u_1u_7)=i_5$, then
$u_5u_4u_3v_1u_6u_7u_1u_2$ is a heterochromatic path of length
$7$. Then we can conclude that there exist vertices
$v_2,v_3\notin\{u_1,u_2,u_3,u_4,u_5,u_6,u_7\}$ such that
$C(u_7v_2)=i_3$, $C(u_7v_3)=i_5$. If $C(u_1u_5)=i_8$, then
$v_1u_3u_2u_1u_5u_6u_7v_2$ is a heterochromatic path of length
$7$, and so we assume $C(u_1u_5)=i_7$. Since
$v_1u_6u_5u_1u_2u_3u_4$ is a heterochromatic path of length $6$,
we have $C(u_1v_1,u_2v_1,u_4v_1,$ $
u_5v_1)-\{i_1,i_2,i_3,i_5,i_6,i_7,i_8\}\neq \emptyset$. Let
\[
P'=\left\{
\begin{array}{ll}
u_4u_3u_2u_1v_1u_6u_7v_3 &\mbox{if $C(u_1v_1)\notin\{i_1,i_2,i_3,i_5,i_6,i_7,i_8\}$;}\\
u_4u_3u_2v_1u_7u_6u_5u_1 &\mbox{if $C(u_2v_1)\notin\{i_1,i_2,i_3,i_5,i_6,i_7,i_8\}$;}\\
v_1u_4u_3u_2u_1u_5u_6u_7 &\mbox{if $C(u_4v_1)\notin\{i_1,i_2,i_3,i_5,i_6,i_7,i_8\}$;}\\
u_4u_3u_2u_1u_5v_1u_6u_7 &\mbox{if
$C(u_5v_1)\notin\{i_1,i_2,i_3,i_5,i_6,i_7,i_8\}$.}
\end{array}
\right.
\]
Then, $P'$ is a heterochromatic path of length $7$, a contradiction.\\
\indent {\bf Case 2.2.1.2} $x_1=3$ and $x_2=6$. Then
$|C(u_1u_7,u_2u_7,u_3u_7,u_4u_7,u5u_7)-\{i_1,i_3,i_4,i_6\}|\geq
4$.\\
\indent If $C(u_1u_7)=i_2$, then $v_1u_3u_4u_5u_6u_7u_1u_2$ is a
heterochromatic path of length $7$; if $C(u_1u_7)=i_5$, then
$u_5u_4u_3v_1u_6u_7u_1u_2$ is a heterochromatic path of length
$7$, a contradiction. So we conclude that
$|C(u_2u_7,u_3u_7,u_4u_7,u_5u_7)-\{i_1,i_3,i_4,i_6\}|=4$. If
$C(u_2u_7)\notin\{i_1,i_3,i_4,i_5,i_6,i_7\}$, then
$v_1u_3u_4u_5u_6u_7u_2u_1$ is a heterochromatic path of length
$7$; if $C(u_2u_7)=i_5$, then $u_5u_4u_3v_1u_6u_7u_2u_1$ is a
heterochromatic path of length $7$, a contradiction. So
$C(u_2u_7)=i_7$. Since $u_1u_2u_3v_1u_6u_5u_4$ is a
heterochromatic path of length $6$, we have
$C(u_1u_3,u_1u_4,u_1u_5,u_1u_6,u_1v_1)-\{i_1,i_2,i_3,i_4,i_5,i_7,i_8\}\neq\emptyset$.
Let
\[
P'=\left\{
\begin{array}{ll}
u_1u_3u_2u_7v_1u_6u_5u_4 &\mbox{if $C(u_1u_3)\notin\{i_1,i_2,i_3,i_4,i_5,i_7,i_8\}$;}\\
u_1u_4u_3u_2u_7v_1u_6u_5 &\mbox{if $C(u_1u_4)\notin\{i_1,i_2,i_3,i_4,i_5,i_7,i_8\}$;}\\
u_1u_5u_4u_3u_2u_7v_1u_6 &\mbox{if $C(u_1u_5)\notin\{i_1,i_2,i_3,i_4,i_5,i_7,i_8\}$;}\\
u_1u_6u_5u_4u_3u_2u_7v_1 &\mbox{if $C(u_1u_6)\notin\{i_1,i_2,i_3,i_4,i_5,i_7,i_8\}$;}\\
u_1v_1u_7u_2u_3u_4u_5u_6 &\mbox{if
$C(u_1v_1)\notin\{i_1,i_2,i_3,i_4,i_5,i_7,i_8\}$.}
\end{array}
\right.
\]
Then, $P'$ is a heterochromatic path of length $7$, a contradiction.\\
\indent {\bf Case 2.2.1.3} $x_1=4$ and $x_2=6$. Then
$|C(u_1u_7,u_2u_7,u_3u_7,u_4u_7,u5u_7)-\{i_1,i_2,i_4,i_6\}|\geq 4$.\\
\indent If $C(u_1u_7)=i_3$, then $v_1u_3u_2u_1u_7u_6u_5u_4$ is a
heterochromatic path of length $7$; if $C(u_1u_7)=i_5$, then
$v_1u_6u_7u_1u_2u_3u_4u_5$ is a heterochromatic path of length
$7$, a contradiction. So we get that
$|C(u_2u_7,u_3u_7,u_4u_7,u_5u_7)-\{i_1,i_2,i_4,i_6\}|=4$. If
$C(u_1u_4)\neq i_7$, then $v_1u_3u_2u_1u_4u_5u_6u_7$ is a
heterochromatic path of length $7$, and so $C(u_1u_4)=i_7$ and
$C(u_1u_6)\notin\{i_1,i_2,i_3,i_4,i_5,i_6,i_7\}$. If
$C(u_1u_7)=i_2$, then $v_1u_3u_4u_5u_6u_7u_1u_2$ is a
heterochromatic path of length $7$, and so there exists a
$v_2\notin\{u_1,u_2,u_3,u_4,u_5,u_6,u_7\}$ such that
$C(u_7v_2)=i_2$. Then, $u_1u_6u_5u_4u_3v_1u_7v_2$ is a
heterochromatic path of length $7$, a contradiction.\\
\indent So, in the case $k=8$, there exists a heterochromatic path
of length $7$ in $G$.

\indent {\bf Case 2.2.2} $k=11$. Denote $i_9=C(u_{y_1}v_1)$,
$i_{10}=C(u_{y_2}v_1)$ and $i_{11}=C(u_{y_3}v_1)$. We distinguish
the following two cases:\\
\indent {\bf Case 2.2.2.1} $y_1=3$, $y_2=6$ and $y_3=8$.\\
\indent We can easily get that $x_3=7$ or $x_3=8$. Since
$3k-3l-2=l-1$ in this case, we have
$|C(u_1u_9,u_2u_9,\ldots,u_6u_9,u_7u_9)-(\{i_1,i_2,i_3,i_4,i_5,i_8\}-\{i_{x_1-1},i_{x_2-1}\})|=7$.
Then, $C(u_1u_9)\in\{i_{x_1-1},i_{x_2-1},i_6,i_7\}$. Let
\[
P'=\left\{
\begin{array}{ll}
v_1u_3u_4u_5u_6u_7u_8u_9u_1u_2 &\mbox{if $C(u_1u_9)=i_2$;}\\
u_4u_5u_6u_7u_8u_9u_1u_2u_3v_1 &\mbox{if $C(u_1u_9)=i_3$;}\\
v_1u_6u_7u_8u_9u_1u_2u_3u_4u_5 &\mbox{if $C(u_1u_9)=i_5$;}\\
u_7u_8u_9u_1u_2u_3u_4u_5u_6v_1 &\mbox{if $C(u_1u_9)=i_6$;}\\
v_1u_8u_9u_1u_2u_3u_4u_5u_6u_7 &\mbox{if $C(u_1u_9)=i_7$.}
\end{array}
\right.
\]
Then, $P'$ is a heterochromatic path of length $9$, a
contradiction. So $C(u_1u_9)=i_4$, and then we can conclude that
$5\in\{x_1,x_2\}$ and $4\notin\{x_1,x_2\}$. Therefore, there
exists a $v_2\notin\{u_1,u_2,\ldots,u_8,u_9\}$ such that
$C(u_9v_2)=i_3$, and $u_5u_6u_7u_8v_1u_3u_2u_1u_9v_2$ is a
heterochromatic path of length $9$, a contradiction.\\
\indent {\bf Case 2.2.2.2} $y_1=3$, $y_2=5$ and $y_3=8$.\\
\indent Since $3k-3l-2=l-1$ in this case, we have
$|\{i_1,i_2,i_5,i_6,i_7,i_8\}\cap\{i_{x_1-1},i_{x_2-1},\\
i_{x_3-1}\}|=2$ and
$|C(u_1u_9,u_2u_9,\ldots,u_6u_9,u_7u_9)-(\{i_1,i_2,i_5,i_6,i_7,i_8\}
-\{i_{x_1-1},i_{x_2-1},\\ i_{x_3-1}\})|=7$. Then we can get that
$x_1=3,x_2=5$ and $x_3=7$, or $x_1=3,x_2=5$ and $x_3=8$, or
$x_1=4,x_2=6$ and $x_3=8$, and
$C(u_1u_9)\in\{i_3,i_4,i_{x_1-1},i_{x_2-1},i_{x_3-1}\}$. Let
\[
P'=\left\{
\begin{array}{ll}
v_1u_3u_4u_5u_6u_7u_8u_9u_1u_2 &\mbox{if $C(u_1u_9)=i_2$;}\\
u_4u_5u_6u_7u_8u_9u_1u_2u_3v_1 &\mbox{if $C(u_1u_9)=i_3$;}\\
v_1u_5u_6u_7u_8u_9u_1u_2u_3u_4 &\mbox{if $C(u_1u_9)=i_4$;}\\
u_6u_7u_8u_9u_1u_2u_3u_4u_5v_1 &\mbox{if $C(u_1u_9)=i_5$;}\\
v_1u_8u_9u_1u_2u_3u_4u_5u_6u_7 &\mbox{if $C(u_1u_9)=i_7$.}
\end{array}
\right.
\]
Then, $P'$ is a heterochromatic path of length $9$, a
contradiction. So $C(u_1u_9)=i_6$, and then $x_1=3,x_2=5,x_3=7$.
Therefore, there exists a $v_2\notin\{u_1,u_2,\ldots,u_8,u_9\}$
such that $C(u_9v_2)=i_5$, and $u_7u_8v_1u_5u_4u_3u_2u_1u_9v_2$ is
a heterochromatic path of length $9$, a contradiction.\\
\indent So, in the case $k=11$, there exists a heterochromatic
path of length $9$ in $G$.

Up to now, we can conclude that if $d^c(v)\geq k\geq 7$ for any
$v\in V(G)$, then $G$ has a heterochromatic path of length at
least $\lceil\frac{2k}{3}\rceil+1$ in $G$. \qed

\section{Long heterochromatic paths under the color neighborhood union condition}

Let $G$ be an edge-colored graph and $s$ a positive integer.
Suppose that $|CN(u)\cup CN(v)|\geq s$ for every pair of vertices
$u$ and $v$ of $G$. It is easy to see that if $s=1,2$ then $G$ has
a heterochromatic path of length $s$, and if $s=3$ then $G$ has a
heterochromatic path of length $2$. In \cite{B-L}, the authors
showed that $G$ has a heterochromatic path of length at least
$\lceil\frac{s}{3}\rceil+1$ for $s>1$. In this section we will
improve this lower bound for $s\geq 4$.

\begin{theo} Let $G$ be an edge-colored graph and $s$ a positive integer.
Suppose that $|CN(u)\cup CN(v)|\geq s\geq 4$ for every pair of
vertices $u$ and $v$ of $G$. Then $G$ has a heterochromatic path
of length at least $\lfloor\frac{2s+4}{5}\rfloor$.
\end{theo}
{\bf\pf} By contradiction. Suppose $P=u_1u_2\ldots u_lu_{l+1}$ is
a longest heterochromatic path of length
$l<\lfloor\frac{2s+4}{5}\rfloor$. Denote $i_j=C(u_ju_{j+1})$
for $j=1,2,\ldots,l$.\\
\indent Since $P$ is a longest heterochromatic path in $G$, there
exist $x_i$'s and $y_j$'s such that $3\leq
x_1<x_2<\ldots<x_{t_1}\leq l+1$ and $2\leq
y_1<y_2<\ldots<y_{t_2}\leq l-1$, and
$t_1=|CN(u_1)-C(P)|=|C(u_1u_{x_1},u_1u_{x_2},\ldots,u_1u_{x_{t_1}})|$,
$t_2=|CN(u_{l+1})-C(P)|=|C(u_{y_1}u_{l+1},u_{y_2}u_{l+1},\ldots,u_{y_{t_2}}u_{l+1})|$
and $C(u_1u_{x_1},u_1u_{x_2},\ldots,u_1u_{x_{t_1}})\cap
C(u_{y_1}u_{l+1},\\ u_{y_2}u_{l+1},\ldots,u_{y_{t_2}}u_{l+1})
=\emptyset$. Then $t_1+t_2\geq s-l>s-\lfloor\frac{2s+4}{5}\rfloor
= \lceil\frac{3s-4}{5}\rceil\geq \frac{3s-4}{5}>\frac{2s-1}{5}\geq
\lfloor\frac{2s+4}{5}\rfloor -1 >l-1$. Denote
$\{z_1,z_2,\ldots,z_{t_3}\}=\{y_1,y_2,\ldots,y_{t_2}\}\cap\{x_1-1,x_2-1,\ldots,x_{t_1}-1\}$,
and so $2\leq z_1<z_2<\ldots<z_{t_3}\leq l-1$. Since $2\leq
y_1<y_2<\ldots<y_{t_2}\leq l-1$ and $2\leq
x_1-1<x_2-1<\ldots<x_{t_1}-1\leq l$, we have $t_3\geq
t_1+t_2-(l-1)>0$. Then, from Lemma 2.2 we can get that
$CN(u_1)-C(u_1u_3,u_1u_4,\ldots,u_1u_l,u_1u_{l+1}) \subseteq
C(P)-\{i_{z_1},i_{z_2},\ldots,i_{z_3}\}$,
$CN(u_{l+1})-C(u_1u_{l+1},u_2u_{l+1},\ldots,u_{l-1}u_{l+1})\subseteq
C(P)-\{i_{z_1},i_{z_2},\ldots,i_{z_3}\}$. So, $CN(u_1)\cup
CN(u_{l+1})$ $\subseteq
(C(P)-\{i_{z_1},i_{z_2},\ldots,i_{z_3}\})\cup C(u_1u_3,u_1u_4,
\ldots,u_1u_l,u_1u_{l+1}, u_2u_{l+1},\ldots,u_{l-1}u_{l+1})$.
Therefore, $|CN(u_1)\cup CN(u_{l+1})|\leq
|C(P)-\{i_{z_1},i_{z_2},\ldots,i_{z_3}\}|+|C(u_1u_3,u_1u_4,\ldots,$
$ u_1u_l,u_1u_{l+1},u_2u_{l+1},\ldots,u_{l-1}u_{l+1})|
=(l-t_3)+(2l-3)=3l-3-t_3\leq
3l-3-(t_1+t_2)+(l-1)=4l-4-(t_1+t_2)\leq 4l-4-(s-l)=5l-4-s <
5*\lfloor\frac{2s+4}{5}\rfloor-4-s\leq s$, a contradiction.

So, if $|CN(u)\cup CN(v)|\geq s\geq 4$ for every pair of vertices
$u$ and $v$ of $G$, then $G$ has a heterochromatic path of length
at least $\lfloor\frac{2s+4}{5}\rfloor$.\qed

Although we cannot show that the above lower bound is best
possible, the following example shows that the best lower bound
cannot be better than $\lfloor\frac{s}{2}\rfloor+1$. Let $s$ be a
positive integer. If $s$ is even, let $G_s$ be the graph obtained
from the complete graph $K_{\frac{s+4}{2}}$ by deleting an edge;
if $s$ is odd, let $G_s$ be the complete graph
$K_{\frac{s+3}{2}}$. Then, color the edges of $G_s$ by different
colors for any two different edges. So, for any $s\geq 1$ we have
that $|CN(u)\cup CN(v)|\geq s$ for any pair of vertices $u$ and
$v$ in $G$, and any longest heterochromatic path in $G$ is of
length $\lfloor\frac{s}{2}\rfloor+1$. This example shows that the
lower bound in our Theorem 4.1 is not very far away from the best.

\end{document}